\documentclass[12pt,reqno]{amsart}

\usepackage[usenames]{color}
\usepackage{amssymb}
\usepackage{amsmath}
\usepackage{amsthm}
\usepackage{amsfonts}
\usepackage{amscd}
\usepackage{graphicx}
\usepackage{mathrsfs}
 
\usepackage{tikz}
\usepackage[colorlinks=true,
linkcolor=webgreen,
filecolor=webbrown,
citecolor=webgreen]{hyperref}

\definecolor{webgreen}{rgb}{0,.5,0}
\definecolor{webbrown}{rgb}{.6,0,0}

\usepackage{color}
\usepackage{fullpage}
\usepackage{float}

\usepackage{graphics}
\usepackage{latexsym}
\usepackage{epsf}
\usepackage{breakurl}

\allowdisplaybreaks
\begin{document}

	\title{	Partial Skew Motzkin Paths }

	\author[H.~Prodinger]{Helmut Prodinger}
	
	\address{Helmut Prodinger,
		Department of Mathematical Sciences, Stellenbosch University,
		7602 Stellenbosch, South Africa, and
		NITheCS (National Institute for Theoretical and Computational Sciences),
		South Africa}
	\email{hproding@sun.ac.za}

	\date{\today}
	
	\begin{abstract}
Motzkin paths consist of up-steps, down-steps, level-steps, and never go below the $x$-axis. They
return to the $x$-axis at the end. The concept of skew Dyck path \cite{Deutsch-italy} is transferred to skew Motzkin paths,
namely, a left step $(-1,-1)$ is additionally allowed, but the path is not allowed to intersect itself.
The enumeration of these combinatorial objects was known \cite{Qing}; here, using the kernel method,
we extend the results by allowing them to end at a prescribed level $j$. The approach is completely based on generating functions.

Asymptotics of the total number of objects as well as the average height are also given.
 	\end{abstract}
	
	
	\maketitle

	\theoremstyle{plain}
	\newtheorem{theorem}{Theorem}
	\newtheorem{corollary}[theorem]{Corollary}
	\newtheorem{lemma}[theorem]{Lemma}
	\newtheorem{proposition}[theorem]{Proposition}
	
	\theoremstyle{definition}
	\newtheorem{definition}[theorem]{Definition}
	\newtheorem{example}[theorem]{Example}
	\newtheorem{conjecture}[theorem]{Conjecture}
	
	\theoremstyle{remark}
	\newtheorem{remark}[theorem]{Remark}

\section{Introduction}

A Dyck path has up-steps and down-steps of one unit each, and cannot go into negative territory. Usually,
one considers returns to the $x$-axis at the end, but for partial paths, this is not required. A Motzkin path
allows additional flat (horizontal) steps of unit length. A skew path allows  `left' step $(-1,-1)$ as well, but 
the path is not allowed to intersect itself. We prefer `red' steps $(1,-1)$, see our analysis in \cite{skew-paper}.
For Motzkin paths, some analysis as provided in \cite{Qing}. Here, we provide further analysis that allows to consider
partial paths as well, so we don't need to land at the $x$-axis. It uses the kernel method \cite{prodinger-kernel}.
Apart from being not below the $x$-axis, the restrictions are that a left (red) step cannot follow or preceed an up-step.

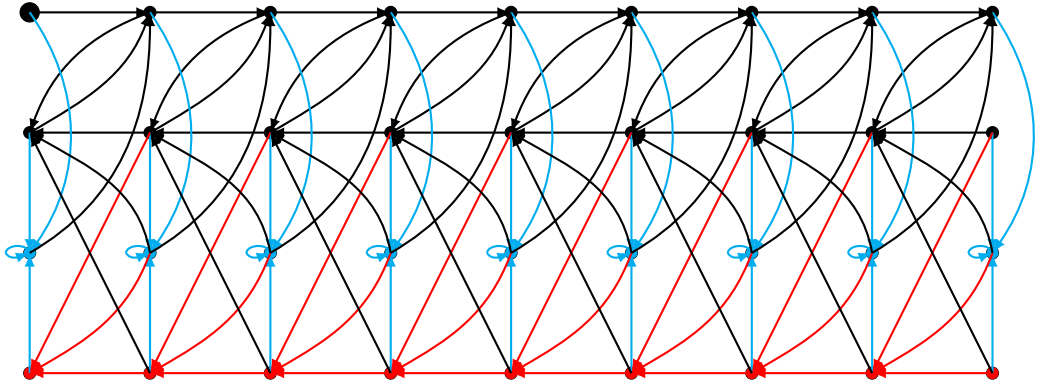
\begin{figure}[h]

	\begin{center}
		\begin{tikzpicture}[scale=1.6]
			\draw (0,0) circle (0.08cm);
			\fill (0,0) circle (0.08cm);
			
			\foreach \x in {0,1,2,3,4,5,6,7,8}
			{
				\draw (\x,0) circle (0.05cm);
				\fill (\x,0) circle (0.05cm);
			}
			
			\foreach \x in {0,1,2,3,4,5,6,7,8}
			{
				\draw (\x,-1) circle (0.05cm);
				\fill (\x,-1) circle (0.05cm);
			}
			
			\foreach \x in {0,1,2,3,4,5,6,7,8}
			{
				\draw (\x,-2) circle (0.05cm);
				\fill[cyan] (\x,-2) circle (0.05cm);
			}
			\foreach \x in {0,1,2,3,4,5,6,7,8}
			{
				\draw (\x,-3) circle (0.05cm);
				\fill[red] (\x,-3) circle (0.05cm);
			}

			\foreach \x in {0,1,2,3,4,5,6,7}
			{
				\draw[ thick,-latex] (\x,0) -- (\x+1,0);
				
			}
			
			\foreach \x in {0,1,2,3,4,5,6,7}
			{
				\draw[ thick,-latex] (\x+1,0) to [out=200, in =70] (\x,-1);
				\draw[ thick,-latex] (\x,-1) to [out=30, in =-110] (\x+1,0);
				\draw[ thick,-latex] (\x+1,-1) to (\x,-1);
			}
		
		\foreach \x in {0,1,2,3,4,5,6,7,8}
		{
			\draw[ thick,cyan,-latex] (\x,0) to  [out=-55, in =55] (\x,-2);%[out=-40, in =40]
\draw[ thick,cyan,-latex] (\x,-1) to   (\x,-2);
\draw[ thick,cyan,-latex] (\x,-3) to   (\x,-2);
		}

		\foreach \x in {0,1,2,3,4,5,6,7}
		{
			\draw[ thick,red,-latex] (\x+1,-1) to  (\x,-3);%[out=-100, in =40] 
				\draw[ thick,red,-latex] (\x+1,-2) to[out=-110, in =30]   (\x,-3);
				\draw[ thick,red,-latex] (\x+1,-3) to   (\x,-3);
		}
	
	\foreach \x in {0,1,2,3,4,5,6,7}
	{
		\draw[ thick,-latex] (\x+1,-3) to    (\x,-1);
	
	}
		\foreach \x in {0,1,2,3,4,5,6,7,8}
	{
		\draw[ thick,cyan ] (\x,-2) to  [out=100, in =80]  (\x-0.2,-2);
		\draw[ thick,cyan,-latex ] (\x-0.2,-2 ) to  [out=-60, in =200] (\x,-2) ;	
	}
	
		\foreach \x in {0,1,2,3,4,5,6,7}
	{
		\draw[ thick,-latex] (\x,-2) to [out=30, in =-90]     (\x+1,-0);
		
	}

\foreach \x in {0,1,2,3,4,5,6,7}
{
	\draw[ thick,-latex] (\x+1,-2) to [out=100, in =-30]     (\x,-1);
	
}

		\end{tikzpicture}
	\end{center}
	\caption{Four layers of states according to the type of steps leading to them. Traditional up-steps and down-steps are black, 
	level-steps are blue, and left steps are red.}
	\label{schoas}
\end{figure}

The situation is best described by a graph (state-diagram); see Figure~\ref{schoas}.

In further sections, the asymptotic equivalent for the number of skew Motzkin paths of given
size is derived, as well as the \emph{height}, meaning that the generating function of paths
with a bounded height (bounded by $H$) is given, as well as the average height, which is
approximately $\textsf{const}\cdot \sqrt{n}$, which is typical for families of paths.
 
 \section{Generating functions for skew Motzkin paths}
 
 We translate the state diagram accordingly; $f_j$, $g_j$, $h_j$, $k_j$ are generating functions in the variable $z$ (marking the length
 of the path), ending at level $j$. The four families are related to the four layers of states. 
 \begin{align*}
f_{j+1}&=zf_j+zg_j+zh_j,\ f_0=1,\\
g_j&=zf_{j+1}+zg_{j+1}+zh_{j+1}+zk_{j+1},\\
h_j&=zf_j+zg_j+zh_j+zk_j,\\
k_j&=zg_{j+1}+zh_{j+1}+zk_{j+1}.
 \end{align*}
Now we introduce bivariate generating functions, namely
\begin{equation*}
F(u):=\sum_{j\ge0}f_ju^j, \ G(u):=\sum_{j\ge0}g_ju^j, \ H(u):=\sum_{j\ge0}h_ju^j, \ K(u):=\sum_{j\ge0}k_ju^j.
\end{equation*}
The recursions then take this form:
 \begin{align*}
	F(u)&=1+zuF(u)+zuG(u)+zuH(u),\\
	uG(u)&=zF(u)+zG(u)+zH(u)+zK(u)-z-zg_0-zh_0-zk_0,\\
	H(u)&=zF(u)+zG(u)+zH(u)+zK(u),\\
	uK(u)&=zG(u)+zH(u)+K(u)-zg_{0}+zh_{0}+zk_{0}.
\end{align*}
Solving the system,
\begin{align*}
	F(u)&=\frac{\mathscr{F}}{2z-u+zu-{z}^{2}u+z{u}^{2}-{z}^{3}-{z}^{3}u} ,\\
G(u)&= \frac{\mathscr{G}}{2z-u+zu-{z}^{2}u+z{u}^{2}-{z}^{3}-{z}^{3}u},\\
H(u)&= \frac{\mathscr{H}}{2z-u+zu-{z}^{2}u+z{u}^{2}-{z}^{3}-{z}^{3}u},\\
K(u)&=\frac{\mathscr{K}}{2z-u+zu-{z}^{2}u+z{u}^{2}-{z}^{3}-{z}^{3}u},
\end{align*}
with $\mathscr{F}=-{z}^{3}+2 z-u+zu+{z}^{2}u+{z}^{2}ug_0+{z}^{2}uh_0+{z}^{2}u
k_0+{z}^{3}ug_0+{z}^{3}uh_0+{z}^{3}uk_0
$, $\mathscr{G}=-{z}^{2}h_0+{z}^{4}+{z}^{4}k_0-{z}^{2}ug_0-{z}^{2}u-{z}
^{2}k_0-{z}^{2}-{z}^{2}uh_0-{z}^{2}uk_0+{z}^{4}h_0
+{z}^{3}-{z}^{2}g_0+zh_0+zg_0+zk_0+{z}^{4}g_0
$, $\mathscr{H}=-{z}^{4}+2 {z}^{2}g_0+2 {z}^{2}h_0+2 {z}^{2}k_0+2 {z}^{2}-zu-{z}^{4}g_0-{z}^{4}h_0-{z}^{4}k_0-{z}^{3}u
g_0-{z}^{3}uh_0-{z}^{3}uk_0$, $\mathscr{K}=zg_0-{z}^{3}-{z}^{2}g_0-{z}^{2}h_0+zk_0-{z}^{2}k_0-{z}^{3}g_0-{z}^{3}
h_0-{z}^{3}k_0$.

One cannot immediately insert $u=0$ to identify the constants, but one can use the kernel method. For that, one factorises the denominator:
\begin{equation*}
2z-u+zu-{z}^{2}u+z{u}^{2}-{z}^{3}-{z}^{3}u=z(u-u_1)(u-u_2)
\end{equation*}
with 
\begin{equation*}
u_1=\frac{1-z+z^2+z^3 +(1+z)W}{2z}, \quad u_2=\frac{1-z+z^2+z^3-(1+z)W}{2z}
\end{equation*}
 and
 \begin{equation*}
W=\sqrt{(1-z)(1-3z-z^2-z^3)}=\sqrt{1-4z+2z^2+z^4}.
 \end{equation*}
Since $u_2\sim 2z$ for small $z$, $u-u_2$ is a `bad' factor and must cancel from both, numerator and denominator. This yields
\begin{align*}
	F(u)&=\frac{-1+z+{z}^{2}+{z}^{2}g_0+{z}^{2}h_0+{z}^{2}k_0+{z}^{3}g_0+{z}^{3}h_0+{z}^{3}k_0	}{z(u-u_1)} ,\\
	G(u)&= \frac{-{z}^{2}-{z}^{2}g_0-{z}^{2}h_0-{z}^{2}k_0}{ z(u-u_1)},\\
	H(u)&= \frac{-z-{z}^{3}g_0-{z}^{3}h_0-{z}^{3}k_0 }{ z(u-u_1)},\\
	K(u)&=\frac{-{z}^{2}g_0-{z}^{2}h_0-{z}^{2}k_0 }{ z(u-u_1)}.
\end{align*}
Now we can plug in $u=0$ and identify the constants:
\begin{align*}
g_0&={\frac {-{z}^{5}+{z}^{3}W-{z}^{2}-zW+3 z-1+W}{2{z}^{2} \left( -2+		{z}^{2} \right) }}
,\\
h_0&=-{\frac {-{z}^{2}+2 z-1+W}{2z}},\\
k_0&=-{\frac {-{z}^{4}+{z}^{3}+{z}^{2}W+zW-3 z+1-W}{2{z}^{2} \left( -2		+{z}^{2} \right) }}.
\end{align*}
Adding these quantities yields
\begin{equation*}
1+g_0+h_0+k_0=-{\frac {-{z}^{2}+2 z-1+W}{2{z}^{2}}},
\end{equation*}
which is the generating function of the number of skew Motzkin paths (returning to the $x$-axis); the series expansion is
\begin{equation*}
1+z+2z^2+5z^3+13z^4+35z^5+97z^6+275z^7+794z^8+2327z^9+6905z^{10}+20705z^{11}+\cdots,
\end{equation*}
as already given in \cite{Qing}, the coefficients are sequence A82582 in \cite{OEIS}.

We further get
\begin{align*}
	F(u)&=\frac{-1+z-z^2-z^3+u_2z}{z(u-u_1)} ,\\
	G(u)&= \frac{(z-u_2)}{ (u-u_1)(1+z)},\\
	H(u)&= \frac{-1-z+2z^2+z^3-zu_2}{ (u-u_1)(1+z)},\\
	K(u)&=\frac{z^2+2z-u_2}{ (u-u_1)(1+z)}.
\end{align*}
Altogether,
\begin{equation*}
F(z)+G(z)+H(z)+K(z)=\frac{-1-z+2z^2+z^3-zu_2}{z(u-u_1)(1+z)}
\end{equation*}
and
\begin{equation*}
[u^j]	\big(F(z)+G(z)+H(z)+K(z)\big)=\frac{1+z-2z^2-z^3+zu_2}{z(1+z)u_1^{j+1}},
\end{equation*}
which is the generating function of partial skew Motzkin paths, landing on level $j$.

Here are the examples for $j=1,2,3,4$ (leading terms only):
\begin{align*}
&z+2z^2+5z^3+13z^4+36z^5+102z^6+295z^7+866z^8+2574z^9+7730z^{10}+23419z^{11},\\
&z^2+3z^3+9z^4+26z^5+77z^6+230z^7+694z^8+2110z^9+6459z^{10}+19890z^{11}+61577z^{12},\\
&z^3+4z^4+14z^5+45z^6+143z^7+451z^8+1421z^9+4478z^{10}+14129z^{11}+44654z^{12},\\
&z^4+5z^5+20z^6+71z^7+242z^8+806z^9+2653z^{10}+8670z^{11}+28213z^{12}.
\end{align*}

One can also substitute $u=1$, which means that \emph{all} partial skew Motzkin paths are counted with respect to length, regardless on which level they end:
\begin{equation*}
 {\frac {2-3 z-7 {z}^{2}-{z}^{3}+{z}^{4}-(2+z) \left( 1+z \right) W}{2z \left( 1+z \right)  \left( 2 {z}^{2}+3
		 z-1 \right) }}.
\end{equation*}
The series expansion is
\begin{equation*}
1+2 z+5 z^2+14 z^3+40 z^4+117 z^5+348 z^6+1049 z^7+3196 z^8+9823 z^9+30413 z^{10}+\cdots
\end{equation*}

\section{Counting flat and left (red) steps}

Using two extra variables $t$ and $w$, we can count the number of flat resp.\ left steps in a skew Motzkin path. The recursions are self-explanatory.
\begin{align*}
	f_{j+1}&=zf_j+zg_j+zh_j,\ f_0=1,\\
	g_j&=zf_{j+1}+zg_{j+1}+zh_{j+1}+zk_{j+1},\\
	h_j&=ztf_j+ztg_j+zth_j+ztk_j,\\
	k_j&=zwg_{j+1}+zwh_{j+1}+zwk_{j+1}.
\end{align*}
Again, here is the system for the multi-variate generating functions;
\begin{align*}
	F(u)&=1+zuF(u)+zuG(u)+zuH(u),\\
	uG(u)&=zF(u)+zG(u)+zH(u)+zK(u)-z-zg_0-zh_0-zk_0,\\
	H(u)&=ztF(u)+ztG(u)+ztH(u)+ztK(u),\\
	uK(u)&=zwG(u)+zwH(u)+zwK(u)-zwg_{0}+zwh_{0}+zwk_{0}.
\end{align*}
And following a similar procedure as before we get
\begin{align*}
	1&+g_0+h_0+k_0=
\frac{-zw+u_2}{z(1+wt)}\\*
&=1+tz+(t^2+1)z^2+(tw+3t+t^3)z^3+(2+6t^2+w+3wt^2+t^4)z^4+\cdots;
\end{align*}
the quantity $u_2$ is now
\begin{equation*}
u_2=\frac{1-tz+w{z}^{2}+tw{z}^{3}-\sqrt{(1-z^2w)(1-2\,tz+ \left( {t}^{2}-4-w \right) {z}^{2}-2tw{z}^{3}-w{t}^{2}{z}^
		{4})}}{2z}.
\end{equation*}
Quantities like $F(u)$, $G(u)$, $H(u)$, $K(u)$ can also be computed easily, following the approach from the previous section.

\section{Asymptotics for the number of skew Motzkin paths}
We must analyze the generating function
\begin{equation*}
	\mathscr{SM}=\frac {(1-z)^2-\sqrt{(1-z)(1-3z-z^2-z^3)}}{2{z}^{2}}
\end{equation*}
which is of the sqrt-type \cite{FO} around the singularity $\rho$ closest to the origin,
which we call  $\rho$. It is a solution of $1-3z-z^2-z^3=0$ and an algebraic number of order 3. We
confine ourselves to real numbers to keep the expressions shorter. 
We expand everything around $\rho=0.295597742522084770980996$.
We find
\begin{align*}
\mathscr{SM}&\sim 	\frac {(1-\rho)^2-\sqrt{2.7142940417132890604(\rho-z)}}{2{\rho}^{2}}\\
&\sim 	2.8392867552141611323-9.4274931376001571585\sqrt{\rho-z}\\
&\sim	2.8392867552141611323- 5.1256244361431546460\sqrt{1-z/\rho}.
\end{align*}
Singularity analysis of generating function \cite{FO} gives the estimate
\begin{equation*}
[z^n]\mathscr{SM}\sim 	5.1256244361431546460\frac{1}{2\sqrt{\pi}}\rho^{-n}n^{-3/2}.
\end{equation*}
The error at $n=100$ is about $3\%$. This is to be expected by this type of approximation.

\section{Skew Motzkin paths of bounded height}

Now we introduce a parameter $H$ and don't allow the path to reach any level higher than $H$. We can still work with the system
 \begin{align*}
	f_{j+1}&=zf_j+zg_j+zh_j,\ 0\le j \le H-1,\ f_0=1,\\
	g_j&=zf_{j+1}+zg_{j+1}+zh_{j+1}+zk_{j+1},\ 0\le j<H, \\
	h_j&=zf_j+zg_j+zh_j+zk_j,\ 0\le j\le H,\\
	k_j&=zg_{j+1}+zh_{j+1}+zk_{j+1},\ 0\le j<H.
\end{align*}
This is now a finite linear system, and we are only interested in paths that return to the $x$-axis. For a given $H$, we write
$s[H]=f_0+g_0+h_0+k_0$ and let Maple compute these quantities for the first 20 values of $H$.

Both, numerator and denominator of $s[H]$ satisfy the recursion
\begin{equation*}
X_{n+2}+(-1+z-z^2-z^3)X_{n+1}+(2z^2-z^4)X_n=0.
\end{equation*}
Thus, adjusting this to the initial conditions, we get
\begin{equation*}
s[n]=\frac{A_o(1+z^3+z^2-z+\omega)^n+ B_o(1+z^3+z^2-z-\omega)^n}
{A_u(1+z^3+z^2-z+\omega)^n+ B_u(1+z^3+z^2-z-\omega)^n}
\end{equation*}
with
\begin{align*}
A_o&=(z^3+z^2+3z-1)(z+1)+(z-1)\omega,\\
B_o&=(z^3+z^2+3z-1)(z+1)-(z-1)\omega,\\
A_u&=(1-z^2)(z^3+z^2+3z-1)+\frac{z^3-z^2+3z-1}{1-z}\omega,\\
B_u&= (1-z^2)(z^3+z^2+3z-1)-\frac{z^3-z^2+3z-1}{1-z}\omega,\\
\omega&=\sqrt{z^6+2z^5+3z^4-5z^2-2z+1}=(1+z)W.
\end{align*}
When $n$ goes to infinity, the second terms go away, and we are left with
\begin{equation*}
s[\infty]=\frac{A_o}{A_u}={\frac {(1-z)^2-\sqrt{(1-z)(1-3z-z^2-z^3)}}{2{z}^{2}}}=\mathscr{SM},
\end{equation*}
as expected.
Now we consider $s[>n]$, the generating function of skew Motzkin paths of height $>n$. Taking differences, we find
\begin{align*}
	s[>n]=s[\infty]-s[n]&=\frac{A_oB_u-A_uB_o}{A_u}\frac{(1+z^3+z^2-z-\omega)^n}
	{A_u(1+z^3+z^2-z+\omega)^n+ B_u(1+z^3+z^2-z-\omega)^n}\\
	&\sim\frac{A_oB_u-A_uB_o}{A^2_u}\frac{\bigg(\dfrac{1+z^3+z^2-z-\omega}{1+z^3+z^2-z+\omega}\bigg)^n}
	{1-\bigg(\dfrac{1+z^3+z^2-z-\omega}{1+z^3+z^2-z+\omega}\bigg)^n}.
	\end{align*}
A computer computation leads to (always in the neighbourhood of $z=\rho$)
\begin{equation*}
\frac{A_oB_u-A_uB_o}{A^2_u}\sim 18.854986275200314363\sqrt{\rho-z}.
\end{equation*}
Now we approximate:
\begin{align*}
	\dfrac{1+z^3+z^2-z-\omega}{1+z^3+z^2-z+\omega}&\sim \frac{0.81760902991166091601-2.1345121404980002137\sqrt{\rho-z}}
	{0.81760902991166091601+2.1345121404980002137\sqrt{\rho-z}}\\
	&\sim \frac{1-2.6106758394395728799\sqrt{\rho-z}}
	{1+2.6106758394395728799\sqrt{\rho-z}}\\
	&\sim 1-5.2213516788791457598\sqrt{\rho-z}\\
	&\sim \exp\bigl({-5.2213516788791457598\sqrt{\rho-z}}\,\bigr)=e^{-t},
\end{align*}
for convenience.
For the average height, we need apart from the leading factor,
\begin{equation*}
\sum_{h\ge0}s[>h] \sim\sum_{h\ge0}\frac{e^{-th}}{1-e^{-th}}.
\end{equation*}
Since we only compute the leading term of the asymptotics of the average height, we might start the sum at $h\ge1$, and expand the geometric series:
\begin{equation*}
\sum_{h\ge1}s[>h] \sim\sum_{h,k\ge1}e^{-thk}=\sum_{k\ge1}d(k)e^{-kt}\sim -\frac{\log t}{t},
\end{equation*}
with $d(k)$ being the number of divisors of $k$. This type of analysis, although having been done often before,
has been described in much detail in \cite{steffel}. Together with the factor in front, we are at
\begin{align*}
&-18.854986275200314363\sqrt{\rho-z}\frac{\log \sqrt{\rho-z}}{5.2213516788791457598\sqrt{\rho-z}}\\
&=	-18.854986275200314363\frac{\log \sqrt{\rho-z}}{5.2213516788791457598}\\
&=-1.8055656307800996608\log(\rho-z)\\
&\sim-1.8055656307800996608\log(1-z/\rho).
\end{align*}
Singularity analysis \cite{FO} gives the following estimate for the coefficient of $z^n$:
\begin{equation*}
	1.8055656307800996608 \frac{\rho^{-n}}{n}.
\end{equation*}
For the average height we need to normalize, which is to divide by the total number of skew Motzkin numbers of size $n$:
\begin{equation*}
\frac{	1.8055656307800996608 \frac{\rho^{-n}}{n}}{5.1256244361431546460\frac{1}{2\sqrt{\pi}}\rho^{-n}n^{-3/2}}
=0.70452513767814089508\sqrt{\pi n}.
\end{equation*}

%\clearpage

\end{document}